\documentclass[runningheads]{llncs}

\usepackage[utf8]{inputenc}
\usepackage{amsmath,amssymb}
\usepackage{graphicx}
\usepackage{booktabs}
\usepackage{hyperref}
\usepackage{enumitem}
\usepackage{geometry}
\usepackage{tikz}
\usetikzlibrary{arrows.meta, positioning, shapes.geometric, fit, backgrounds}
\usepackage{pgfplots}
\pgfplotsset{compat=1.17}
\usepackage{tabularx}
\usepackage{array}
\usepackage{adjustbox}
\begin{document}

\title{When Protection Hurts: A Cytokine Storm Theory of Supply Chain Overreaction}
\author{Gurkirat Wadhwa}
\authorrunning{Gurkirat et al.}
\institute{Industrial Engineering and Operations Research (IEOR), IIT Bombay\\
\email{gurkirat@iitb.ac.in}}

\maketitle

%=========================================================================%

\begin{abstract}
Supply chain resilience research has traditionally viewed defensive actions such as safety stock, backup suppliers, expedited transportation, and emergency sourcing as beneficial. However, disruptions such as the COVID-19 pandemic, the 2021 semiconductor shortage, and the Suez Canal blockage highlight how synchronized hoarding, duplicate ordering, and capacity reservation can amplify congestion and operational costs, sometimes worsening the original disruption.
Motivated by this observation, we develop a decision-theoretic framework in which Protective Response Intensity (PRI) is an endogenous decision variable. We introduce Response-Induced Loss (RIL) to capture costs generated by defensive actions and characterize a unique optimal response that balances disruption mitigation against self-induced losses. We derive a closed-form critical response threshold beyond which RIL exceeds the original disruption intensity, a phenomenon termed a *Supply Chain Cytokine Storm*. We further propose the Cytokine Storm Index (CSI), measuring response-induced loss relative to the original disruption, and discuss its potential calibration using operational data.
Extending the framework to interconnected supply networks, we derive a spectral stability condition distinguishing bounded decentralized responses from network-wide propagation. Numerical experiments show that optimal responses remain below the storm threshold, whereas proportional response heuristics can exceed it. Network responses converge or diverge according to the derived stability condition. Overall, our findings highlight that supply chain resilience depends not only on mitigating disruptions but also on keeping defensive responses proportionate to the shocks they are intended to address.
\end{abstract}

\textbf{Keywords:} Supply chain resilience; behavioral overreaction; response-induced loss; network cascades; convex optimization.

%=========================================================================%

\section{Introduction}
\label{sec:intro}

Supply chain resilience is traditionally viewed through a simple principle: greater preparedness leads to lower disruption losses. Consequently, firms invest in safety stock, backup suppliers, flexible capacity, expedited transportation, and emergency sourcing to improve their ability to withstand unexpected disruptions \cite{tang2006,christopher2004,sheffi2005}. Implicit in this perspective is the assumption that defensive actions are inherently beneficial and that increasing preparedness can only improve resilience.

Recent large-scale disruptions suggest that this assumption is incomplete. During the COVID-19 pandemic, the global semiconductor shortage, and the Suez Canal blockage, firms responded to uncertainty with synchronized hoarding, duplicate ordering, capacity reservation, and aggressive expediting. Although these actions were individually rational, their collective effect amplified congestion, prolonged shortages, and increased operational costs well beyond those caused by the original disruption itself \cite{shih2020,handfield2020,attinasi2021,tran2025suez,notteboom2021}. These episodes indicate that supply chains may suffer not only from external shocks, but also from excessive defensive responses.

This phenomenon closely resembles a \emph{cytokine storm} in immunology, where the immune system's response to a pathogen becomes so aggressive that it causes more damage than the infection itself \cite{tisoncik2012,mehta2020}. Motivated by this analogy, we investigate an analogous question for supply chains: \emph{Can protective actions themselves become the dominant source of operational loss?}. This analogy is used for motivating the intuition behind the proposed framework.

To address this question, we model a firm's \emph{Protective Response Intensity} (PRI) as a decision variable rather than a fixed parameter and explicitly account for the costs generated by defensive actions through the notion of \emph{Response-Induced Loss} (RIL). Building on this formulation, we introduce the \emph{Cytokine Storm Index} (CSI), which quantifies the proportion of loss attributable to the response itself, derive a closed-form threshold distinguishing beneficial defensive actions from excessive ones, and extend the framework to interconnected supply networks by establishing a spectral condition governing whether defensive responses remain localized or propagate throughout the network.

Our work differs from existing resilience and behavioral supply chain literature in three important respects. First, unlike traditional resilience models that treat preparedness as uniformly beneficial \cite{tang2006,christopher2004,sheffi2005}, we explicitly model protective actions as optimization decisions that balance mitigation benefits against self-induced operational costs. Second, whereas the bullwhip-effect literature demonstrates how ordering variability propagates through supply chains \cite{lee1997,sterman1989,croson2006}, we derive a closed-form analytical condition identifying when defensive responses themselves become counterproductive. Third, by extending the analysis to interconnected supply networks, we show that excessive responses are not solely firm-level phenomena but are fundamentally shaped by network interactions, linking resilience analysis with spectral methods commonly used in network science \cite{wang2003,pastor2015}.

The main contributions of this paper are as follows. We formulate a decision-theoretic framework that models Protective Response Intensity (PRI) as an endogenous decision variable and introduces Response-Induced Loss (RIL) to explicitly capture the operational costs generated by defensive actions. Building on this framework, we establish the existence and uniqueness of the optimal response, derive its comparative statics, characterize a closed-form threshold separating beneficial protection from excessive defensive behavior, and propose the Cytokine Storm Index (CSI) as a quantitative measure of self-induced operational loss. Finally, we extend the analysis to interconnected supply networks by deriving a spectral stability condition governing the propagation of defensive responses and demonstrate the resulting threshold and cascade phenomena through numerical experiments.

The remainder of the paper is organized as follows. Section~\ref{sec:related} reviews the related literature. Section~\ref{sec:model} presents the proposed model, Section~\ref{sec:analysis} develops the theoretical analysis, Section~\ref{sec:network} extends the framework to interconnected supply networks, Section~\ref{sec:numerical} presents the numerical experiments, and Section~\ref{sec:conclusion} concludes the paper.

%=========================================================================%
\section{Related Work}
\label{sec:related}

Resilience research typically models a firm's or network's ability to absorb, adapt to, and recover from disruptions, with defensive investment entering as a parameter that improves absorption or shortens recovery time \cite{tang2006,christopher2004,sheffi2005}. Tang \cite{tang2006} catalogs mitigation strategies---postponement, flexible supply base, economic incentives---as levers firms should pull as far as budget allows; Christopher and Peck \cite{christopher2004} and Sheffi \cite{sheffi2005} similarly frame resilience as a capacity to be built up, with redundancy and flexibility appearing as unambiguously protective. These models describe a single firm acting against a well-defined risk reasonably well, but they share an implicit assumption that resilience-enhancing investment is weakly beneficial, leaving no room for the response itself to become a source of loss.

A second stream shows that overreaction can emerge even without any individual decision-maker behaving irrationally. The bullwhip effect literature shows that order variance amplifies upstream, driven by demand-signal processing, order batching, price fluctuation, and rationing, even under individually optimal ordering policies \cite{lee1997}. Laboratory studies using the beer distribution game show that human subjects systematically over-order relative to the optimal benchmark, especially just after a shock \cite{sterman1989,croson2006}. This literature establishes overreaction as real and measurable, but characterizes it indirectly, through order-variance decompositions under a fixed ordering policy, rather than by asking when a \emph{chosen} level of responsiveness becomes counterproductive; no closed-form threshold analogous to ours appears in this line of work. Nor does this literature offer a metric comparable to CSI: standard resilience diagnostics such as fill rate, time-to-recovery, or bullwhip order-variance ratios summarize how well a system absorbs or reacts to a shock, but to the best of our knowledge, these studies do not isolate the fraction of loss that the response itself creates. CSI is intended to fill exactly that gap, and comparing it empirically against these established metrics is a natural next step we flag in Section~\ref{sec:numerical}.

A third stream draws on immune-inspired methods and network cascades. Artificial Immune Systems have long served as metaheuristics for supplier selection, scheduling, and routing \cite{dasgupta2011}, and terms like ``organizational immunity'' occasionally appear in the practitioner literature as informal metaphors for adaptive capacity, stopping short of formalizing an \emph{overactive} response. Separately, network science and epidemiology use spectral conditions to determine whether a linear process on a network converges or diverges: a threshold of the form $\beta\rho(A)\lessgtr 1$, with $\rho(A)$ the spectral radius of the network's adjacency matrix, governs whether a contagion dies out or spreads \cite{wang2003,pastor2015}. We adapt this machinery to a different substrate: what propagates in our network model is a firm's own defensive reaction, not a pathogen or a default, and it is the unchecked growth of that reaction that constitutes the harm.

Taken together, these streams motivate the need for a unified decision-theoretic perspective. While existing studies have substantially advanced our understanding of disruption propagation, resilience strategies, and behavioral amplification, relatively less attention has been devoted to characterizing when protective actions themselves become counterproductive. To the best of our knowledge, we are not aware of an analytical framework that explicitly identifies a critical threshold between beneficial protection and self-induced operational loss, while also examining how such behavior propagates over supply networks. Section~\ref{sec:model} develops a framework to address these questions.

%=========================================================================%
\section{Model}
\label{sec:model}

\subsection{Setup}

\begin{figure}[h]
\centering
\begin{tikzpicture}[
node distance=0.9cm,
>=Stealth,
box/.style={
draw,
rounded corners=2mm,
minimum width=3.7cm,
minimum height=0.9cm,
align=center,
font=\small,
fill=blue!8
},
greenbox/.style={
draw,
rounded corners=2mm,
minimum width=3.8cm,
minimum height=0.9cm,
align=center,
font=\small,
fill=green!15
},
redbox/.style={
draw,
rounded corners=2mm,
minimum width=3.8cm,
minimum height=0.9cm,
align=center,
font=\small,
fill=red!15
},
graybox/.style={
draw,
rounded corners=2mm,
minimum width=4.2cm,
minimum height=0.9cm,
align=center,
font=\small,
fill=gray!12
},
arrow/.style={->,thick}
]

%%%%%%%%%%%%%%%%%%%%%%%%%%%%%%%%%%%%%%

\node[box] (D)
{External Disruption\\$D$};

\node[box,below=of D]
(R)
{Protective Response Intensity\\$R$};

%%%%%%%%%%%%%%%%%%%%%%%%%%%%%%%%%%%%%%

\node[greenbox,
below left=1.4cm and 1.8cm of R]
(LD)
{Mitigated Exogenous Loss\\
$L_D(R)=De^{-\mu R}$};

\node[redbox,
below right=1.4cm and 1.8cm of R]
(LR)
{Response-Induced Loss\\
$L_R(R)=\theta R^2$};

%%%%%%%%%%%%%%%%%%%%%%%%%%%%%%%%%%%%%%

\node[graybox,
below=2.3cm of R]
(L)
{Total Loss\\
$L(R)=L_D(R)+L_R(R)$};

%%%%%%%%%%%%%%%%%%%%%%%%%%%%%%%%%%%%%%

\node[box,
below=of L]
(CSI)
{Cytokine Storm Index\\
$\displaystyle
CSI=\frac{L_R(R)}{D}
$};

%%%%%%%%%%%%%%%%%%%%%%%%%%%%%%%%%%%%%%

\node[greenbox,
below left=1.5cm and 1.8cm of CSI]
(Safe)
{\textbf{Resilient Region}\\
$CSI\le1$};

\node[redbox,
below right=1.5cm and 1.8cm of CSI]
(Storm)
{\textbf{Supply Chain}\\
\textbf{Cytokine Storm}\\
$CSI>1$};

%%%%%%%%%%%%%%%%%%%%%%%%%%%%%%%%%%%%%%

\node[graybox,
below=2.5cm of CSI]
(Network)
{Network Response\\
$\beta\rho(A)<1$ Stable \qquad
$\beta\rho(A)\ge1$ Cascading};

%%%%%%%%%%%%%%%%%%%%%%%%%%%%%%%%%%%%%%

\draw[arrow] (D)--(R);

\draw[arrow] (R)--(LD);

\draw[arrow] (R)--(LR);

\draw[arrow] (LD)--(L);

\draw[arrow] (LR)--(L);

\draw[arrow] (L)--(CSI);

\draw[arrow] (CSI)--(Safe);

\draw[arrow] (CSI)--(Storm);

\draw[arrow] (Storm)--(Network);

\draw[arrow] (Safe)--(Network);
\end{tikzpicture}

\caption{Model Description}
\label{fig:causal-chain}
\end{figure}
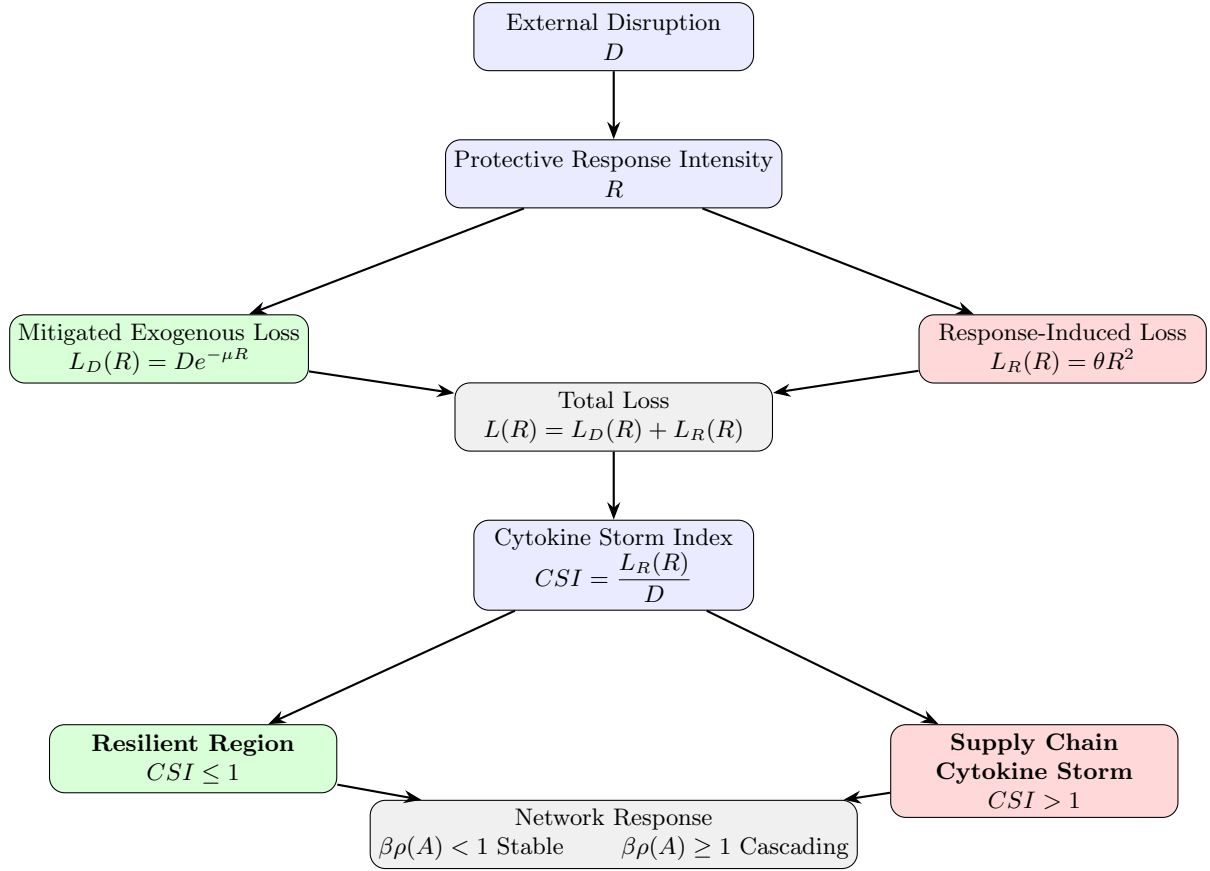

Consider a firm facing a disruption of intensity $D \geq 0$, representing, for example, a normalized capacity loss, demand shock, or lead-time increase expressed on a common cost-equivalent scale. To mitigate the disruption, the firm selects a \emph{Protective Response Intensity} (PRI), denoted by $R \geq 0$, which represents the aggregate level of defensive actions such as increasing safety stock, expediting shipments, reserving transportation capacity, and engaging emergency suppliers. While these actions reduce the losses directly attributable to the disruption, they also incur operational and coordination costs. Consequently, the firm's response creates a trade-off between disruption mitigation and response-induced costs. Figure~\ref{fig:causal-chain} illustrates the framework describing the model.

\begin{definition}[Mitigated exogenous loss]
The exogenous loss actually realized after mitigation is
\begin{equation}
L_D(R) = D\, e^{-\mu R},
\label{eq:mitigated-loss}
\end{equation}
where $\mu > 0$ is the \emph{mitigation effectiveness} parameter. A larger response reduces exogenous loss, with diminishing returns.
\end{definition}

\begin{definition}[Response-Induced Loss]
The Response-Induced Loss (RIL) is
\begin{equation}
L_R(R) = \theta R^2,
\label{eq:ril}
\end{equation}
where $\theta > 0$ is the \emph{overreaction cost coefficient}, capturing hoarding, order amplification, and congestion costs that grow faster than linearly in the response.
\end{definition}
The quadratic response cost in \eqref{eq:ril} captures increasing marginal costs of defensive actions, including congestion, coordination costs, and capacity shortages. Together, \eqref{eq:mitigated-loss} and \eqref{eq:ril} capture the trade-off between diminishing mitigation benefits and increasing overreaction costs while preserving analytical tractability.
\begin{definition}[Total loss and Cytokine Storm Index]
Total loss is
\begin{equation}
L(R) = L_D(R) + L_R(R) = D e^{-\mu R} + \theta R^2,
\label{eq:total-loss}
\end{equation}
and the \emph{Cytokine Storm Index} (CSI) at response level $R$ is
\begin{equation}
\mathrm{CSI}(R) = \frac{L_R(R)}{D} = \frac{\theta R^2}{D}.
\label{eq:csi}
\end{equation}
\end{definition}
The Cytokine Storm Index (CSI) measures the magnitude of self-inflicted loss relative to the original disruption intensity $D$, rather than relative to the mitigated (realized) exogenous loss $L_D(R)$. This normalization is deliberate: $D$ is fixed when the response is chosen, whereas $L_D(R)$ already reflects that response. Thus CSI provides a response-independent benchmark of proportionality: $\mathrm{CSI}(R)<1$ means response-induced loss is below the original disruption intensity, while $\mathrm{CSI}(R)>1$ indicates disproportionate response-induced loss. CSI is therefore a benchmark-relative measure of response excessiveness rather than a decomposition of realized total loss. 
\begin{definition}[Supply Chain Cytokine Storm]
\label{def:storm}
A firm is in a \emph{Supply Chain Cytokine Storm} at response level $R$ if
\begin{equation}
\mathrm{CSI}(R) > 1 \quad \Longleftrightarrow \quad \theta R^2 > D.
\label{eq:storm-condition}
\end{equation}
\end{definition}
A Supply Chain Cytokine Storm occurs when the response-induced loss exceeds the original disruption intensity $D$ that motivated the response, not necessarily the (already-mitigated) exogenous loss actually realized once the response is in place; the two coincide only when $\mu=0$. Moreover, this does not imply that the response is ineffective---on the contrary, the response may substantially mitigate exogenous losses while simultaneously creating even greater operational losses through excessive defensive actions. Thus, the storm condition captures disproportionate response intensity relative to the size of the original disruption, rather than inadequate disruption mitigation.
\subsection{Optimal response}

Given the loss function in \eqref{eq:total-loss}, a rational firm seeks to determine the level of Protective Response Intensity that minimizes its total expected loss. The resulting optimization balances the mitigation benefits of defensive actions against the operational costs they generate, thereby identifying the response that achieves the best trade-off between preparedness and overreaction. To highlight the value of this optimal policy, we compare it with a simple behavioral benchmark in which firms respond proportionally to the observed disruption, i.e., $R=\alpha D$ for a fixed responsiveness parameter $\alpha>0$. This heuristic represents panic-driven or overly aggressive ordering behavior frequently associated with the bullwhip effect. Section~\ref{sec:analysis} characterizes the optimal response analytically and establishes the conditions under which both decision rules enter the cytokine storm region.
%=========================================================================%
\section{ Analysis}
\label{sec:analysis}

We now analyze the proposed model and establish its fundamental theoretical properties. We first prove the existence and uniqueness of the optimal response, followed by its comparative statics with respect to the model parameters. We then characterize the critical threshold separating beneficial and excessive defensive responses, and finally compare the behavior of the optimal policy with that of a simple proportional response heuristic.
\begin{proposition}[Convexity and unique optimal response]
\label{prop:convex}
The loss function $L(R) = D e^{-\mu R} + \theta R^2$ is strictly convex on $R \geq 0$ for any $D,\mu,\theta>0$, and therefore has a unique minimizer $R^\star \geq 0$, characterized by
\begin{equation}
2\theta R^\star = \mu D\, e^{-\mu R^\star}.
\label{eq:foc}
\end{equation}
\end{proposition}

\noindent\textit{Proof.}
Differentiating, $L'(R) = -\mu D e^{-\mu R} + 2\theta R$ and $L''(R) = \mu^2 D e^{-\mu R} + 2\theta$. Since $D,\mu,\theta>0$, $L''(R)>0$ everywhere, so $L$ is strictly convex. At $R=0$, $L'(0)=-\mu D<0$, and $L'(R)\to+\infty$ as $R\to\infty$ since the linear term dominates the vanishing exponential; by the intermediate value theorem there exists $R^\star>0$ with $L'(R^\star)=0$, i.e.\ \eqref{eq:foc}. Strict convexity implies this is the unique global minimizer.
\qed
Proposition~\ref{prop:convex} establishes that neither insufficient nor excessive defensive actions are optimal. Instead, the firm faces a well-defined trade-off between disruption mitigation and response-induced costs, resulting in a unique optimal level of protection. Equation~\eqref{eq:foc} implicitly defines $R^\star$ as a function of $D,\mu,\theta$; it can also be written in closed form via the Lambert $W$ function as $R^\star = \tfrac{1}{\mu}W\!\big(\tfrac{\mu^2 D}{2\theta}\big)$, but we work directly with the implicit condition \eqref{eq:foc} throughout, since it is the more transparent form for the comparative statics that follow. It is the threshold $R_c$ in Proposition~\ref{prop:threshold} below, not $R^\star$ itself, that we refer to elsewhere as ``closed-form.''
\begin{proposition}[Comparative statics]
\label{prop:comparative-statics}
The optimal response $R^\star(D,\theta)$ implicitly defined by \eqref{eq:foc} satisfies $\partial R^\star/\partial D > 0$ and $\partial R^\star/\partial \theta < 0$.
\end{proposition}

\noindent\textit{Proof.}
Write $F(R,D,\theta) = 2\theta R - \mu D e^{-\mu R} = 0$. By Proposition~\ref{prop:convex}, $F_R = 2\theta+\mu^2De^{-\mu R} = L''(R)>0$ at $R=R^\star$. By the implicit function theorem, $\partial R^\star/\partial D = \mu e^{-\mu R^\star}/F_R>0$ and $\partial R^\star/\partial\theta = -2R^\star/F_R<0$.
\qed
Proposition~\ref{prop:comparative-statics} provides an intuitive characterization of the optimal response. Firms facing more severe disruptions should adopt stronger defensive actions, whereas higher overreaction costs warrant greater restraint. Thus, the optimal response by the firms is jointly determined by disruption severity and the economic consequences of defensive actions.
\begin{proposition}[Critical response threshold]
\label{prop:threshold}
For fixed $D,\theta>0$, define $R_c = \sqrt{D/\theta}$. Then $\mathrm{CSI}(R)>1$ iff $R>R_c$; a Supply Chain Cytokine Storm occurs exactly when the response exceeds $R_c$.
\end{proposition}

\noindent\textit{Proof.} Immediate from Definition~\ref{def:storm}: $\theta R^2>D \iff R>\sqrt{D/\theta}=R_c$. \qed

\begin{corollary}[Storm-free optimum]
\label{cor:no-storm}
For every $D,\mu,\theta>0$, the optimal response satisfies $R^\star<R_c$. Equivalently,
\[
\mathrm{CSI}(R^\star)<1.
\]
Thus, under the proposed loss function, an optimally-controlled firm never enters a Supply Chain Cytokine Storm.
\end{corollary}

\noindent\textit{Proof.}
By Proposition~\ref{prop:convex}, $L'(R)$ is strictly increasing and
$L'(R^\star)=0$. It therefore suffices to show that $L'(R_c)>0$.
Using $R_c=\sqrt{D/\theta}$,
\[
L'(R_c)
=
-\mu D e^{-\mu R_c}+2\theta R_c
=
2\sqrt{D\theta}
-\mu D e^{-\mu\sqrt{D/\theta}}.
\]
Let
\[
x=\mu\sqrt{\frac{D}{\theta}}>0.
\]
Then
\[
L'(R_c)
=
\sqrt{D\theta}\left(2-xe^{-x}\right).
\]
Since $xe^{-x}\leq 1/e<2$ for every $x>0$, we have
$L'(R_c)>0$. Because $L'$ is strictly increasing and
$L'(R^\star)=0$, it follows that $R^\star<R_c$. Hence
$\mathrm{CSI}(R^\star)<1$.
\qed

The corollary shows that cytokine storms cannot arise from the unconstrained optimal response under the proposed loss function. They instead require departures from optimal decision-making, such as the proportional response heuristic analyzed below, or additional constraints and network interactions that amplify defensive responses.

\begin{proposition}[Storm threshold under proportional response]
\label{prop:heuristic-threshold}
If a firm follows $R=\alpha D$ for fixed $\alpha>0$, a storm occurs iff
\begin{equation}
\alpha > \alpha_c := \frac{1}{\sqrt{D\theta}}.
\label{eq:alpha-c}
\end{equation}
\end{proposition}

\noindent\textit{Proof.} Substituting $R=\alpha D$ into Proposition~\ref{prop:threshold}: $\alpha D>\sqrt{D/\theta} \iff \alpha>1/\sqrt{D\theta}$. \qed

For $D=10$, $\theta=0.1$, this gives $\alpha_c=1$: any responsiveness at or above one-for-one with the observed disruption triggers a storm regardless of $\mu$, since $\mu$ never enters \eqref{eq:alpha-c}. Section~\ref{sec:numerical} illustrates this and every other proposition above numerically. We next extend the single-firm threshold to a network of interacting firms, since none of the three motivating episodes in Section~\ref{sec:intro} involved a single isolated decision-maker.

%=========================================================================%

\section{Network Extension}
\label{sec:network}

The analysis in the previous sections considers a single decision maker that chooses its protective response independently. In real supply chains, however, firms rarely operate in isolation. A firm's defensive actions---such as inventory hoarding, emergency sourcing, or aggressive capacity reservation can alter the operating environment faced by its suppliers and customers, often inducing them to increase their own defensive responses. Consequently, overreaction is not merely an individual optimization problem but also a network phenomenon in which local decisions may propagate throughout the supply chain.

To capture this mechanism, we extend the proposed framework from a single firm to a network of interacting firms. Rather than developing a full multi-agent game in which firms jointly optimize their responses, which will be a part of future work, we isolate the propagation mechanism through which protective responses spread across supply chain relationships. The network model is therefore a reduced-form behavioral propagation model rather than a strategic equilibrium. The complete multi-agent game analysis of multiple firms is an interesting future direction.

Suppose the supply network consists of $n$ firms connected through a weighted interaction graph. Firm $i$ experiences a disruption of intensity $D_i$ and incurs the same response trade-off introduced in the single-firm model. Aggregating these individual losses yields the network objective
\begin{equation}
L(\mathbf{R})
=
\sum_{i=1}^{n}
\left(
D_i e^{-\mu R_i}
+
\theta R_i^2
\right),
\label{eq:network-loss}
\end{equation}
which is the natural multi-firm extension of the single-firm loss function in \eqref{eq:total-loss}. Likewise, we define the network-level Cytokine Storm Index as
\begin{equation}
\mathrm{CSI}(\mathbf{R})
=
\frac{\sum_{i=1}^{n}\theta R_i^2}
{\sum_{i=1}^{n}D_i},
\label{eq:network-csi}
\end{equation}
which measures the aggregate proportion of response-induced loss across the entire supply network.

Let $\mathcal N=\{1,\ldots,n\}$ denote the set of firms and let $A\in\mathbb{R}_{+}^{n\times n}$ be a weighted, nonnegative, irreducible adjacency matrix whose element $A_{ij}$ quantifies the influence of firm $j$'s response on firm $i$. Each firm first determines a baseline response
\[
b_i=\alpha D_i,
\]
representing the proportional heuristic introduced in Proposition~\ref{prop:heuristic-threshold}. This initial response is subsequently amplified through interactions with neighboring firms according to
\begin{equation}
R_i
=
\alpha D_i
+
\beta
\sum_j
A_{ij}R_j,
\label{eq:network-response}
\end{equation}
where $\beta\ge0$ is a behavioral contagion parameter describing the strength with which firms imitate or react to the defensive actions of their neighbors. In matrix notation,
\[
\mathbf R
=
\mathbf b
+
\beta A\mathbf R.
\]
The following theorem characterizes when these behavioral interactions converge to a finite equilibrium and when they instead generate cascading defensive responses throughout the supply network.

\begin{theorem}[Network cascade condition]
\label{thm:cascade}
Let $\rho(A)$ be the spectral radius of the irreducible, nonnegative
matrix $A$. If $\beta\rho(A)<1$, the system has a unique nonnegative
solution
\begin{equation}
R^\star
=
(I-\beta A)^{-1}b
=
\sum_{k=0}^{\infty}\beta^k A^k b,
\label{eq:neumann}
\end{equation}
and the dynamics
\begin{equation}
R^{(t+1)}=b+\beta A R^{(t)},
\qquad
R^{(0)}=0,
\end{equation}
converge to $R^\star$. If $\beta\rho(A)\geq1$, the same dynamics
diverge whenever $b\geq0$ and $b\neq0$.
\end{theorem}

\noindent\textit{Proof.}
Suppose first that $\beta\rho(A)<1$. Then
\[
\rho(\beta A)=\beta\rho(A)<1.
\]
Hence, by the Neumann series theorem,
\[
(I-\beta A)^{-1}
=
\sum_{k=0}^{\infty}(\beta A)^k
=
\sum_{k=0}^{\infty}\beta^k A^k
\]
converges, and therefore $I-\beta A$ is invertible. Since
$A\geq0$ and $b\geq0$, every term in the series is nonnegative.
Consequently,
\[
R^\star
=
\sum_{k=0}^{\infty}\beta^k A^k b
\geq0.
\]
Moreover, starting from $R^{(0)}=0$, repeated substitution in the
dynamics gives
\[
R^{(t)}
=
\sum_{k=0}^{t-1}(\beta A)^k b
=
\sum_{k=0}^{t-1}\beta^k A^k b.
\]
Thus $R^{(t)}$ is precisely the $t$-th partial sum of the convergent
Neumann series, and hence
\[
R^{(t)}\longrightarrow R^\star
\qquad\text{as }t\to\infty.
\]

To establish uniqueness, let $\widetilde{R}$ be any solution of
\[
\widetilde{R}=b+\beta A\widetilde{R}.
\]
Then
\[
(I-\beta A)\widetilde{R}=b.
\]
Since $I-\beta A$ is invertible,
\[
\widetilde{R}
=
(I-\beta A)^{-1}b
=
R^\star.
\]
Hence the solution is unique, and it is nonnegative.

Now suppose that $\beta\rho(A)\geq1$. Since $A$ is nonnegative and
irreducible, the Perron--Frobenius theorem guarantees the existence
of a strictly positive left Perron eigenvector $w>0$ satisfying
\[
w^\top A=\rho(A)w^\top.
\]
Because $b\geq0$ and $b\neq0$, and because every component of $w$
is strictly positive,
\[
w^\top b>0.
\]

Multiplying the recursion
\[
R^{(t+1)}=b+\beta A R^{(t)}
\]
from the left by $w^\top$ gives
\[
w^\top R^{(t+1)}
=
w^\top b+\beta w^\top A R^{(t)}.
\]
Using $w^\top A=\rho(A)w^\top$, we obtain
\[
w^\top R^{(t+1)}
=
w^\top b+\beta\rho(A)w^\top R^{(t)}.
\]
Since $R^{(0)}=0$, iteration yields
\begin{equation}
w^\top R^{(t)}
=
(w^\top b)
\sum_{k=0}^{t-1}
\bigl(\beta\rho(A)\bigr)^k.
\label{eq:perron_growth}
\end{equation}

If $\beta\rho(A)>1$, the geometric sum in
\eqref{eq:perron_growth} diverges geometrically. If
$\beta\rho(A)=1$, then
\[
w^\top R^{(t)}
=
t\,w^\top b
\longrightarrow\infty.
\]
Thus, in either case,
\[
w^\top R^{(t)}\longrightarrow\infty.
\]
Since $w>0$ and $R^{(t)}\geq0$, the divergence of
$w^\top R^{(t)}$ implies that the response vector $R^{(t)}$ is
unbounded. Therefore, the network response diverges whenever
$\beta\rho(A)\geq1$.
\qed

\newcommand{\ignore}[1]{}
\ignore{
\begin{theorem}[Network cascade condition]
\label{thm:cascade}
Let $\rho(A)$ be the spectral radius of the irreducible, nonnegative matrix $A$. If $\beta\rho(A)<1$, the system has a unique nonnegative solution
\begin{equation}
R^\star = (I-\beta A)^{-1} b = \sum_{k=0}^\infty \beta^k A^k b,
\label{eq:neumann}
\end{equation}
and the dynamics $R^{(t+1)}=b+\beta A R^{(t)}$, $R^{(0)}=0$, converge to $R^\star$. If $\beta\rho(A)\geq1$, the same dynamics diverge whenever $b$ has a nonzero component along the Perron eigenvector of $A$.
\end{theorem}

\noindent\textit{Proof.}
If $\beta\rho(A)<1$, all eigenvalues of $\beta A$ lie inside the unit disc, so $I-\beta A$ is invertible and $(I-\beta A)^{-1}=\sum_k(\beta A)^k$ converges. Since $A,b\geq0$, every partial sum is nonnegative, and $R^{(t)}=\sum_{k=0}^{t-1}(\beta A)^k b$ is exactly the $t$-th partial sum, so $R^{(t)}\to R^\star$. If $\beta\rho(A)\geq1$, let $v>0$ be the Perron eigenvector for $\lambda_{\max}=\rho(A)$, guaranteed unique and strictly positive by irreducibility, and write $b=cv+b_\perp$, $c\neq0$; the component of $R^{(t)}$ along $v$ grows as $c\sum_{k=0}^{t-1}(\beta\lambda_{\max})^k$, unbounded as $t\to\infty$.
\qed}

Theorem~\ref{thm:cascade} shows that whether decentralized overreaction dies out or explodes is governed by a single scalar condition, structurally identical to the epidemic thresholds reviewed in Section~\ref{sec:related}, except that what spreads here is a defensive reaction, not a pathogen. As the network-level loss in \eqref{eq:network-loss} aggregates $\theta R_i^2$ across firms, an unbounded $R^{(t)}$ corresponds to an unbounded, network-wide storm even when every individual $D_i$ is small: a network-wide storm need not require a network-wide shock, only a sufficiently connected, reactive network responding to a shock concentrated in one firm---precisely the dynamic reported in the semiconductor and Suez episodes, where a shock local to one link of the chain produced congestion that spread well beyond it.

%=========================================================================%
\section{Numerical Illustration and Discussion}
\label{sec:numerical}

The theoretical results established in Sections~\ref{sec:analysis}--\ref{sec:network} characterize the conditions under which defensive responses remain beneficial, become counterproductive, and propagate across interconnected firms. This section illustrates these analytical insights through a sequence of numerical experiments. Specifically, we examine whether an optimally controlled firm remains below the cytokine storm threshold (Corollary~\ref{cor:no-storm}), whether a proportional response heuristic crosses the analytically derived threshold $\alpha_c=1/\sqrt{D\theta}$ (Proposition~\ref{prop:heuristic-threshold}), and whether network responses converge or diverge according to the spectral condition $\beta\rho(A)\lessgtr1$ (Theorem~\ref{thm:cascade}).
\subsection{Optimal response remains below the storm threshold}

Consider first a firm that follows the optimal response policy. Table~\ref{tab:optimal} reports the corresponding optimal response $R^\star$ and Cytokine Storm Index (CSI) over a range of disruption intensities $D$ and overreaction cost coefficients $\theta$, with the mitigation effectiveness parameter fixed at $\mu=0.30$. The storm-free result of Corollary~\ref{cor:no-storm} applies to every parameter combination in the table. 
\begin{table}[h]
\centering
\label{tab:optimal}
\renewcommand{\arraystretch}{1.15}

\begin{tabularx}{\linewidth}{
>{\centering\arraybackslash}X
>{\centering\arraybackslash}X
>{\centering\arraybackslash}X
>{\centering\arraybackslash}X
>{\centering\arraybackslash}X
>{\centering\arraybackslash}X}
\toprule
$\mathbf{D}$ &
$\boldsymbol{\theta}$ &
$\mathbf{R^\star}$ &
$\mathbf{L^\star}$ &
$\mathbf{CSI(R^\star)}$ &
\textbf{Storm?} \\
\midrule
5  & 0.05 & 4.224 & 2.300 & 0.178 & No \\
10 & 0.05 & 5.597 & 3.432 & 0.157 & No \\
15 & 0.05 & 6.467 & 4.246 & 0.139 & No \\
20 & 0.05 & 7.110 & 4.897 & 0.126 & No \\
\midrule
5  & 0.10 & 3.026 & 2.933 & 0.183 & No \\
10 & 0.10 & 4.224 & 4.600 & 0.178 & No \\
15 & 0.10 & 5.008 & 5.847 & 0.167 & No \\
20 & 0.10 & 5.597 & 6.863 & 0.157 & No \\
\midrule
5  & 0.20 & 2.036 & 3.544 & 0.166 & No \\
10 & 0.20 & 3.026 & 5.865 & 0.183 & No \\
15 & 0.20 & 3.704 & 7.681 & 0.183 & No \\
20 & 0.20 & 4.224 & 9.201 & 0.178 & No \\
\bottomrule
\end{tabularx}
\vspace{2mm}
\caption{Optimal response $R^\star$, optimal loss $L^\star$, and Cytokine Storm Index (CSI) for different disruption intensities and overreaction costs ($\mu=0.30$).}
\end{table}

Consistent with Corollary~\ref{cor:no-storm}, the optimal response remains below the cytokine storm threshold throughout the entire parameter range considered, with $\mathrm{CSI}(R^\star)$ never exceeding $0.19$, far below the critical value of $1$. Moreover, the comparative statics predicted by Proposition~\ref{prop:comparative-statics} are clearly reflected in Table~\ref{tab:optimal}: for a fixed overreaction cost coefficient $\theta$, the optimal response $R^\star$ increases with the disruption intensity $D$, whereas for a fixed disruption intensity $D$, $R^\star$ decreases as the overreaction cost coefficient $\theta$ increases.
 Corollary~\ref{cor:no-storm} establishes the result globally, the numerical results illustrate  the storm-free property of the optimal response. We next analyze what happens when firms replace this optimal decision rule with a proportional heuristic.
\subsection{Heuristic response can storm}

\begin{figure}[h!]
\centering
\begin{tikzpicture}
\begin{axis}[
    width=0.82\textwidth,
    height=4.4cm,
    xlabel={Response intensity $R$},
    ylabel={Loss},
    xmin=0, xmax=20,
    ymin=0, ymax=45,
    legend pos=north west,
    legend style={font=\footnotesize},
    grid=both
]
\addplot[domain=0:20, samples=100, thick, blue] {10*exp(-0.30*x)};
\addplot[domain=0:20, samples=100, thick, red] {0.10*x*x};
\addplot[domain=0:20, samples=100, thick, black, dashed] {10*exp(-0.30*x)+0.10*x*x};
\addplot[mark=*, mark size=2.2pt, only marks, black] coordinates {(4.224,4.600)};
\addplot[mark=triangle*, mark size=2.6pt, only marks, teal] coordinates {(10,10)};
\legend{$L_D(R)$, $L_R(R)$, $L(R)$, Optimum $R^\star$, Storm crossing $R_c$}
\end{axis}
\end{tikzpicture}
\caption{Loss decomposition for $D=10$, $\theta=0.10$, $\mu=0.30$. $L_R(R)$ crosses the constant level $D=10$ at $R_c=10$, the storm boundary of Proposition~\ref{prop:threshold}. $R^\star\approx4.22$ lies far to the left of this boundary.}
\label{fig:loss-curve}
\end{figure}
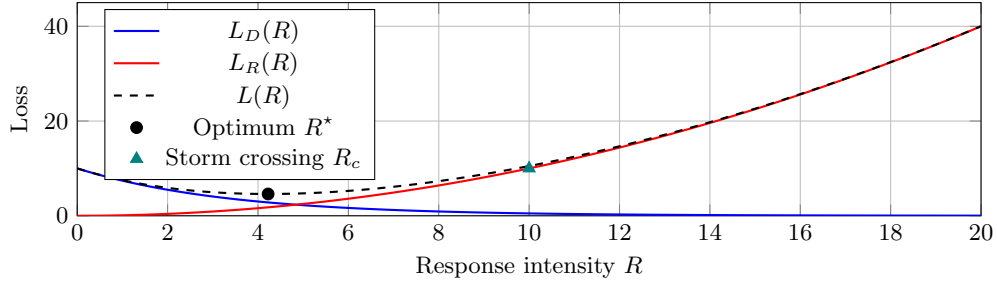

Now suppose the firm abandons that trade-off and reacts in proportion to the shock it sees. We fix   $D=10$, $\theta=0.10$, so that $\alpha_c=1/\sqrt{D\theta}=1.0$ , and vary the proportional responsiveness $\alpha$ in $R=\alpha D$ to obtain Table \ref{tab:heuristic}.
\begin{table}[h]
\centering
\label{tab:heuristic}
\renewcommand{\arraystretch}{1.15}

\begin{tabularx}{\linewidth}{
>{\centering\arraybackslash}X
>{\centering\arraybackslash}X
>{\centering\arraybackslash}X
>{\centering\arraybackslash}X
>{\centering\arraybackslash}X}
\toprule
$\boldsymbol{\alpha}$ &
$\mathbf{R=\alpha D}$ &
$\mathbf{L(R)}$ &
$\mathbf{CSI(R)}$ &
\textbf{Storm?} \\
\midrule
0.10 & 1.0  & 7.508  & 0.010 & No \\
0.30 & 3.0  & 4.966  & 0.090 & No \\
0.50 & 5.0  & 4.731  & 0.250 & No \\
0.80 & 8.0  & 7.307  & 0.640 & No \\
1.10 & 11.0 & 12.469 & 1.210 & Yes \\
1.40 & 14.0 & 19.750 & 1.960 & Yes \\
1.60 & 16.0 & 25.682 & 2.560 & Yes \\
2.00 & 20.0 & 40.025 & 4.000 & Yes \\
\bottomrule
\end{tabularx}
\vspace{2mm}
\caption{Loss and Cytokine Storm Index (CSI) under the proportional response heuristic $R=\alpha D$ for $D=10$, $\theta=0.10$, and $\mu=0.30$. The analytical storm threshold occurs at $\alpha_c=1.0$.}
\end{table}
Loss is not monotonic in $\alpha$: by Proposition~\ref{prop:convex}, $L(R)$ is strictly convex in $R$ and hence in $\alpha$ for fixed $D$, so it falls to a minimum near $\alpha\approx0.5$---close to the optimal $R^\star=4.224$ from Table~\ref{tab:optimal}, i.e.\ $\alpha^\star\approx0.42$---then rises sharply as overreaction costs dominate. The storm boundary falls exactly at $\alpha=1.0$: at $\alpha=1.10$, $\mathrm{CSI}=1.21$ already exceeds one. Figure~\ref{fig:loss-curve} plots $L(R)$ against $R$ with both components, marking the crossing $R_c=\sqrt{D/\theta}=10$ and the optimum $R^\star$ side by side.

Although responding proportionally to the magnitude of an observed disruption may appear to be a reasonable defensive strategy, the numerical results show that such a heuristic can readily push a firm beyond the cytokine storm threshold. In contrast, the optimal response policy remains below this threshold by explicitly balancing mitigation benefits against response-induced costs. These findings highlight that effective supply chain resilience depends not on maximizing the intensity of defensive actions, but on calibrating them to the underlying disruption while accounting for the costs generated by the response itself.
\subsection{Network cascade simulation}

We next consider the network extension of Section~\ref{sec:network}, where the question is no longer how one firm reacts, but what happens when many firms watch and react to each other. We generate a random undirected network of $n=12$ firms (edge probability $0.15$, symmetrized, connected so that $A$ is irreducible), for which $\rho(A)=3.044$. Firm~1 has baseline reaction $\alpha D_1=5.0$; all others have $\alpha D_i=1.0$. We simulate $R^{(t+1)}=b+\beta A R^{(t)}$, $R^{(0)}=0$, for $t=0,1,\dots,30$, at $\beta\rho(A)=0.8$ (below the critical value $1/\rho(A)=0.329$) and $\beta\rho(A)=1.2$ (above it), so the two runs sit on opposite sides of the threshold in Theorem~\ref{thm:cascade}.

The numerical study illustrates the transition between convergent and cascading defensive responses predicted by Theorem~\ref{thm:cascade}. Robustness across network topologies, sizes, and contagion levels is left for future work.
Table~\ref{tab:cascade} reports the total network response $\sum_i R_i^{(t)}$ at every fifth period from $t=0$ to $t=30$, so that both regimes are compared at the same, evenly spaced set of time points.

\begin{figure}[h!]
\centering
\begin{tikzpicture}
\begin{axis}[
    width=0.8\textwidth,
    height=4.4cm,
    xlabel={Time period $t$},
    ylabel={Total network response},
    ymin=0, ymax=1000,
    legend pos=north west,
    legend style={font=\footnotesize},
    grid=both
]
\addplot[mark=*, thick, blue, mark size=1.3pt] coordinates
{(0,0.00) (1,16.00) (2,25.99) (3,34.27) (4,40.66) (5,45.80) (6,49.85)
 (7,53.11) (8,55.69) (9,57.77) (10,59.43) (11,60.75) (12,61.81)
 (13,62.66) (14,63.34) (15,63.88) (16,64.32) (17,64.66) (18,64.94)
 (19,65.16) (20,65.34) (21,65.48) (22,65.60) (23,65.69) (24,65.76)
 (25,65.82) (26,65.87) (27,65.90) (28,65.93) (29,65.96) (30,65.98)};
\addplot[mark=triangle*, thick, red, mark size=1.3pt] coordinates
{(0,0.00) (1,16.00) (2,30.98) (3,49.63) (4,71.19) (5,97.17) (6,127.97)
 (7,165.05) (8,209.27) (9,262.49) (10,326.11) (11,402.66) (12,494.27)
 (13,604.42) (14,736.36) (15,894.94)};
\legend{$\beta\rho(A)=0.8$ (converges to $66.05$), $\beta\rho(A)=1.2$ (diverges, shown to $t{=}15$)}
\end{axis}
\end{tikzpicture}
\caption{Total network response over time for the twelve-firm network. Below threshold, the response flattens out and converges to the closed-form limit $66.05$; above it, it grows without bound (shown to $t=15$ to keep the plot's $y$-axis readable --- by $t=30$ this trajectory reaches over $14{,}000$).}
\label{fig:cascade}
\end{figure}
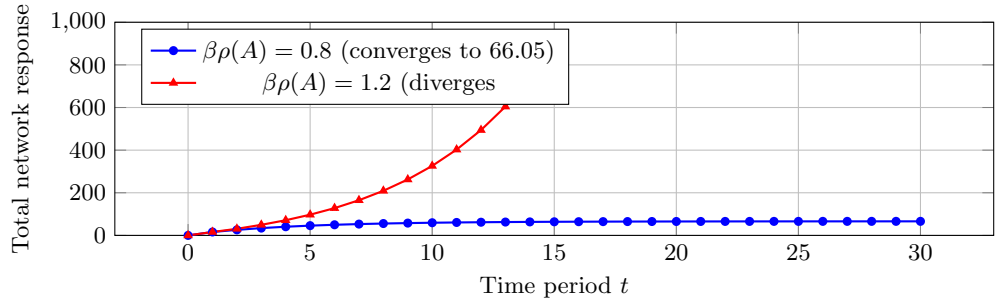

\begin{table}[h]
\centering
\label{tab:cascade}
\renewcommand{\arraystretch}{1.15}

\begin{tabular}{ccc}
\toprule
\textbf{Time} $t$
&
\textbf{Convergent Regime}
($\beta\rho(A)=0.8$)
&
\textbf{Divergent Regime}
($\beta\rho(A)=1.2$)
\\
\midrule
0  & 0.00  & 0.00     \\
5  & 45.80 & 97.17    \\
10 & 59.43 & 326.11   \\
15 & 63.88 & 894.94   \\
20 & 65.34 & 2{,}309.75 \\
25 & 65.82 & 5{,}830.75 \\
30 & 65.98 & 14{,}591.48 \\
\bottomrule
\end{tabular}
\vspace{2mm}
\caption{Evolution of the total network response}
\end{table}
Below threshold, the total network response converges toward the closed-form limit $(I-\beta A)^{-1}b$ from \eqref{eq:neumann}, which evaluates to $66.05$ for this network; by $t=30$ the simulated total of $65.98$ is within $0.1\%$ of this limit, and the remaining gap keeps shrinking geometrically at rate $\beta\rho(A)=0.8$ per period. Above threshold, the response instead grows by a roughly constant multiplicative factor each period ($\approx\beta\rho(A)=1.2$, as the per-period ratio in the simulation converges to this value from above), matching the dominant-eigenvalue argument in the proof of Theorem~\ref{thm:cascade}: by $t=30$ it has reached over $14{,}000$, about $16\times$ its $t=15$ value. Figure~\ref{fig:cascade} plots both trajectories at every period, making the divergence between the two regimes visible from around $t\approx 8$ onward.

In interconnected supply chains, resilience depends not only on how individual firms respond to disruptions but also on how those responses propagate through the network. Even a well-intentioned local reaction can escalate into a system-wide overreaction when amplified by highly connected partners. Consequently, managers should promote firm-level resilience initiatives with mechanisms that limit behavioral contagion, such as information sharing, coordinated decision-making, and greater supply chain visibility.

Taken together, the numerical results are consistent with the theoretical analysis. The optimal response remains below the cytokine storm threshold across the tested parameter range, whereas simple proportional response heuristics readily enter the storm region. At the network level, defensive responses either remain contained or propagate into system-wide overreaction depending on the spectral stability condition of Theorem~\ref{thm:cascade}. These findings suggest that disruption impacts can arise not only from the external shock itself but also from the collective response of firms. This observation is consistent with events following the Suez Canal blockage, where widespread inventory building, rerouting, and capacity reservation contributed to congestion and prolonged recovery across global supply chains \cite{notteboom2021}.

\subsection{Managerial Implications}

The framework yields four managerial implications: (i)  the firms should optimize rather than maximize defensive actions, since marginal response costs eventually dominate mitigation benefits; (ii) monitor $R_c=\sqrt{D/\theta}$ and CSI as warning measures of disproportionate response; (iii) evaluate resilience using both disruption-induced and response-induced losses; and (iv) coordinate responses across supply-network partners, since behavioral contagion can amplify individually reasonable actions into network-wide overreaction.
\subsection{Limitations}

There are several limitations. First, the reduced-form parameters $\theta$, $\mu$, and $\beta$ summarize overreaction, mitigation effectiveness, and behavioral contagion but are not estimated from industry-specific data; estimating them from historical disruption events is an important next step. Second, CSI is normalized using the original disruption intensity $D$ rather than the mitigated exogenous loss, providing a response-independent benchmark; alternative normalizations and comparisons with fill rate, recovery time, and bullwhip measures would further clarify its practical value. Third, the network extension is a reduced-form behavioral propagation model rather than a strategic equilibrium, and its numerical illustration uses a single representative network realization rather than a systematic sweep across topologies. Finally, the network extension assumes static topology and exogenous disruption characteristics; dynamic supply networks with jointly optimized defensive responses is an interesting direction for future research.
%=========================================================================%

\section{Conclusion}
\label{sec:conclusion}

We develop a decision-theoretic framework for analyzing when protective actions in supply chains become counterproductive. By introducing the concepts of \emph{Protective Response Intensity} (PRI), \emph{Response-Induced Loss} (RIL), and the \emph{Cytokine Storm Index} (CSI), the proposed framework explicitly distinguishes losses caused by external disruptions from those generated by excessive defensive responses. The analysis establishes a unique optimal response, derives a closed-form threshold separating proportionate and excessive defensive actions, and extends the framework to interconnected supply networks through a spectral stability condition governing the propagation of behavioral responses.
The numerical experiments illustrate the theoretical results. The optimal response remains below the cytokine storm threshold, whereas simple proportional response heuristics readily enter the storm region. The network experiments further illustrate how behavioral responses either remain localized or propagate into system-wide overreaction depending on the underlying network interactions.
Overall, the proposed framework suggests that effective supply chain resilience depends not only on mitigating external disruptions but also on regulating the intensity of defensive responses.

\bibliographystyle{splncs04}

\end{document}